# Uncovering Ramanujan's "Lost" Notebook: An Oral History

ROBERT P. SCHNEIDER

ABSTRACT. Here we weave together interviews conducted by the author with three prominent figures in the world of Ramanujan's mathematics, George Andrews, Bruce Berndt and Ken Ono. The article describes Andrews's discovery of the "lost" notebook, Andrews and Berndt's effort of proving and editing Ramanujan's notes, and recent breakthroughs by Ono and others carrying certain important aspects of the Indian mathematician's work into the future. Also presented are historical details related to Ramanujan and his mathematics, perspectives on the impact of his work in contemporary mathematics, and a number of interesting personal anecdotes from Andrews, Berndt and Ono.

*"I still say to myself when I am depressed and find myself forced to listen to pompous and tiresome people, 'Well, I have done one thing you could never have done, and that is to have collaborated with Littlewood and Ramanujan on something like equal terms.'"*

—G. H. Hardy [**12**]

The backdrop is Spring, 1976. Beyond the bookshelves, marble busts and tall library windows, the green manicured lawns of Trinity College are framed with flowers. Stone gargoyles and medieval halls, blackened by time, stand as guardians to bright young students whizzing past on bicycles and lounging on the bank of the river Cam [**11**].

Within the archives of the college's Wren Library, an American mathematician opens a box of papers deposited seven years earlier. Over five decades after these papers were mailed to England from India, their story has only begun; on page after page of mathematics, George Andrews recognizes the handwriting of number theorist Srinivasa Ramanujan [**3**]. He has discovered previously unevaluated work from the final year of Ramanujan's brief life. Thus begins a thrilling chapter of contemporary mathematics that is still being written [**14**].

In the mythos of mathematics, the romantic story of Ramanujan is unparalleled. In 1914 the impoverished, self-taught Hindu genius departed from his home of Kumbakonam, India, to sail to the hallowed Trinity College in Cambridge, England. There he was to collaborate with the preeminent British mathematician of his age, G. H. Hardy [**13**], with whom he had exchanged letters filled with exquisite, nearly unbelievable mathematics. Five years later, as a world-famous researcher and Fellow of the Royal Society, Ramanujan returned to India as a national hero—albeit suffering from malnutrition allegedly brought on by the difficulty of adhering to a strictly vegetarian Hindu diet in post-war Britain. Tragically, he passed away one year later at age 32 from a mysterious illness, now believed to have been a curable intestinal parasitic infection [**17**].

With Andrews's finding of this "lost" notebook, not truly lost but languishing unread for more than 50 years, a flood of new ideas was released into the modern world [**14**]. The notes Andrews discovered had traveled a tangled path leading from the Indian mathematician's young widow Janaki Ammal, who gathered the papers after Ramanujan's death [**6**], through the hands of prominent mathematicians such as Hardy, G. N. Watson and R. A. Rankin, before alighting in the archives of Trinity [**4**].

Whereas Ramanujan's earlier work dealt largely with classical number-theoretic objects such as $q$-series, theta functions, partitions and prime numbers—exotic, startling, breathtaking identities built up from infinite series, integrals and continued fractions—in these newfound papers, Andrews found never-before-seen work on the mysterious "mock theta functions" hinted at in a letter written to Hardy in Ramanujan's final months, pointing to realms at the very edge of the mathematical landscape. The content of Ramanujan's lost notebook is too rich, too ornate, too strange to be developed within the scope of the present article. We provide a handful of stunning examples below, intended only to tantalize—perhaps mystify—the reader, who is encouraged to let his or her eyes wander across the page, picking patterns like spring flowers from the wild field of symbols.

The following are two fantastic $q$-series identities found in the lost notebook, published by Andrews soon after his discovery [**2**], in which $q$ is taken to be a complex number with $|q| < 1$:

$$1 + \sum_{n=1}^{\infty} \frac{q^n}{\prod_{j=1}^{n}(1 + aq^j)\left(1 + \frac{q^j}{a}\right)}$$
$$= (1 + a)\sum_{n=0}^{\infty} a^{3n} q^{\frac{1}{2}n(3n+1)}(1 - a^2 q^{2n+1}) - a\frac{\sum_{n=0}^{\infty}(-1)^n a^{2n} q^{n(n+1)/2}}{\prod_{j=1}^{\infty}(1 + aq^j)\left(1 + \frac{q^j}{a}\right)};$$

$$\left(\cfrac{1}{1 + \cfrac{q}{1 + \cfrac{q^2}{1 + \cfrac{q^3}{1 + \cfrac{q^4}{\vdots}}}}}\right)^3 = \frac{\sum_{n=0}^{\infty} q^{5n^2+4n}\frac{1 + q^{5n+2}}{1 - q^{5n+2}} - \sum_{n=0}^{\infty} q^{5n^2+6n+1}\frac{1 + q^{5n+3}}{1 - q^{5n+3}}}{\sum_{n=0}^{\infty} q^{5n^2+2n}\frac{1 + q^{5n+1}}{1 - q^{5n+1}} - \sum_{n=0}^{\infty} q^{5n^2+8n+3}\frac{1 + q^{5n+4}}{1 - q^{5n+4}}}.$$

Another surprising expression [**14**] involves an example of a mock theta function provided by Ramanujan in the final letter he sent to Hardy

$$f(q) = \sum_{n=0}^{\infty} a_f(n) q^n := \sum_{n=0}^{\infty} \frac{q^{n^2}}{(1 + q)^2(1 + q^2)^2 \cdots (1 + q^n)^2}$$

$$= 1 + q - \cdots + 17503 q^{99} + \cdots.$$

In the words of mathematician Ken Ono, a contemporary trailblazer in the field of mock theta functions, "Obviously Ramanujan knew much more than he revealed [**14**]." Indeed, Ramanujan then "miraculously claimed" that the coefficients of this mock theta function obey the asymptotic relation

$$a_f(n) \sim \frac{(-1)^{n-1}}{2\sqrt{n - \frac{1}{24}}} \cdot e^{\pi\sqrt{\frac{n}{6} - \frac{1}{144}}}.$$

The new realms pointed to by the work of Ramanujan's final year are now understood to be ruled by bizarre mathematical structures known as harmonic Maass forms [**14**]. This broader perspective was only achieved in the last ten years, and has led to cutting-edge science, ranging from cancer research to the physics of black holes to the completion of group theory [**13**].

Yet details of George Andrews's unearthing of Ramanujan's notes are only sparsely sketched in the literature; one can detect but an outline of the tale surrounding one of the most fruitful mathematical discoveries of our era. In hopes of contributing to a more complete picture of this momentous event and its significance, here we weave together excerpts from interviews we conducted with Andrews [**1**] documenting the memories of his trip to Trinity College, as well as from separate interviews with mathematicians Bruce Berndt [**6**] and Ken Ono [**13**], who have both collaborated with Andrews in proving and extending the contents of Ramanujan's famous lost notebook.

**GEORGE ANDREWS:** I did not know that this was there—this lost notebook—so it came as a complete surprise.

**KEN ONO:** It's a miracle event. The discovery of the lost notebook sparked an explosion in research, not just on the part of George, but on the part of people who had been working in this area.

**BRUCE BERNDT:** It was a startling discovery, and naturally I was extremely curious as to what was there.

**ONO:** The way to think of the lost notebook is—well, what was Ramanujan? He was a great anticipator of theories that would be developed long after his time.

**BERNDT:** Ramanujan is generally regarded as the greatest Indian mathematician in history. "Ramanujan"— you really should get the pronunciation of it in Indian. I think all of us in the West still slightly mispronounce it, because it's very hard for us to say a name without emphasizing a syllable. In the correct Indian pronunciation there is no emphasis on the syllables, so in other words, "Ra-man-u-jan" is more even. I mean, I say "Ra-*man*-u-jan" so that there's emphasis on the second syllable.

**ONO:** I first became aware of Ramanujan at a fairly young age, as a high school student. My father is a number theorist, and I learned of Ramanujan from him. My father told me about the legendary tale of this great untrained mathematician from South India and it almost didn't sound believable, it sounds more like an old fable or an ancient Chinese tale.

**ANDREWS:** Ramanujan is a unique character in the history of mathematics because he emerged from nowhere—actually from poverty in Southern India—and discovered things about mathematics that surprised major mathematicians in England like Hardy and Littlewood.

**ONO:** The story of Ramanujan is one of the truly romantic legends in modern science, the story of an untrained mathematician who discovered the love of mathematics on his own—not because of coursework or because of the media—discovered the beauty of pure mathematics on his own and for whatever reason, reasons we'll never understand, had a gift, a gift for discovering formulas that would go on to anticipate generations of mathematicians well beyond his time.

**BERNDT:** He received an invitation from G. H. Hardy to come to England so that his mathematical talents could be developed.

**ONO:** Ramanujan's home is in the town of Kumbakonam. Kumbakonam is today what I would consider a fairly large city. Certainly when Ramanujan grew up there, it would have been a smaller town. If you were to visit Ramanujan's home today—or in the early twentieth century—what you would find is a small, one- or two-room structure. There is a bedroom that opens up to the main street, Sarangapani Street. It's a very busy, narrow street with lots of cars and people on bicycles and motorcycles. This is chaos—this is modern India, this country with over a billion people. During his time, it was still a fairly busy commercial street—imagine little shops and merchants lining the whole three- or four-block long street.

Behind his room there is a little kitchen area, where of course there is no running water, there is no sink, there's a little hallway that leads to a back courtyard where one will find a well. This is his whole home—very small, very damp, dark when there is no sunlight. During the summer months, I imagine it would be oppressively hot and humid. The front of the house has a little porch, the porch where we are told that Ramanujan did much of his work sitting on the bench, with slate on his lap scrawling away, madly scribbling formulas. When you step off the porch and turn right, you will be facing Ramanujan's temple, Sarangapani Temple.

Hardy invited Ramanujan to study with him in Cambridge, with the idea that he would offer him proper training in the Western ways of doing modern mathematics. So Ramanujan accepted the invitation. At first he declined because of religious reasons—you know, the trip from South India to Cambridge was quite long. By any event, he accepted that invitation.

**BERNDT:** He lived with E. H. Neville for a couple of weeks after arriving, and then some rooms became available at Trinity College, and so he had a room there.

**ANDREWS:** His subsequent collaboration with Hardy was path-breaking work for the number theory of the twentieth century. Hardy had seen clearly from his correspondence with Ramanujan what a surprising, unique genius this man was. They worked for several years.

**ONO:** As we know it, Ramanujan produced over 30 papers [under Hardy] [**15**].

**ANDREWS:** It definitely wasn't a [typical] student-teacher relationship, it was two very different colleagues collaborating. They published many papers. One of their papers is the basis of the study of what is called probabilistic number theory. Another of their papers, giving a formula for the partition function, is the beginning of what is called the circle method which has been so important in analytic number theory throughout the twentieth century. [Ramanujan's] work on the tau function has led to numerous things—the tau function is basically the coefficient in the expansion of a simple infinite product [**8**]. You take a variable $x$ and then you take the infinite product $x \prod_{n=1}^{\infty} (1 - x^n)^{24}$; there are various reasons why this should be a very important and significant function, and it turns out to be such. Also, one could fairly say that the tremendous advances in the theory of modular forms in the twentieth century began with studying the work of Ramanujan and some of the conjectures he left.

**ONO:** Most of the papers were solo papers by Ramanujan, but it was clear that Hardy had a major impact on helping develop Ramanujan into a professional mathematician, someone who could write papers that could be peer-reviewed, and so on. It's also clear that Hardy had an impact on Ramanujan, because Ramanujan chose to work on some problems in number theory which were central to the mathematicians in Cambridge and Europe at the time. So Hardy's influence was clear in terms of helping to develop [Ramanujan's] professionalism, and in suggesting topics of research.

**BERNDT:** I think there are certain things that Ramanujan was concerned with and was working on, that in fact enabled Hardy to expand his [own] interests. I mean, there are various things that Hardy wouldn't have done if it weren't for Ramanujan.

**ONO:** Apart from the papers that Ramanujan wrote, it's very clear that Hardy's mentorship did not do anything but perhaps enhance Ramanujan's creativity. Some of the papers that Ramanujan wrote during his years in Cambridge are the seminal papers for many different subjects, which mathematicians over the next 100 years would develop into big theories.

**BERNDT:** He left India in March of 1914 for England, and stayed there five years before departing to return home, I think in hope that his illness would subside if he returned to a warmer climate and better food. Unfortunately that did not turn out to be the case.

**ONO:** The legacy of those five years in England is definitely alive today.

**BERNDT:** He died slightly more than a year after he returned to India.

**ONO:** He returned to India as a hero. I think it's important to remember that Ramanujan was an Indian, a South Indian who had made it very big in the Western world of science. And as an Indian who had lived his entire life under British imperialistic rule, he'd reached the greatest heights that a scientist could reach in England. He was elected Fellow of the Royal Society as an Indian—that was a big deal. And so the Indians were very proud of him. He was treated and hailed as a hero when he returned to India, and it's rather tragic that he died at this early age of 32.

**BERNDT:** He had been diagnosed with tuberculosis both in England and then India. But about maybe 15 years ago a British physician, D. A. B. Young, looked at all of the records that he could find and all of the descriptions of Ramanujan's illness that were made by Hardy and others, and concluded that Ramanujan died of hepatic amoebiasis, which is an amoebic infection of the liver if dysentery is not treated properly. It would have been curable. Even in the late stages, if the diagnosis had been properly made, Ramanujan could have been saved.

**ONO:** Unfortunately, as many legends go, Ramanujan was discovered as a young man, but he also died as a young man, long before his time. What we as mathematicians really love about this story, is that, well, it's really a story about an artist—Ramanujan was an artist with functions, he was an artist with numbers. He had a way with them that experts in Western Europe hadn't figured out. He was able to overcome—supersede—the accumulated wisdom of centuries of work on the part of mathematicians in England and in Germany. He did this all as an untrained mathematician, as a two-time college dropout. And as legends go, his ideas still continue to inspire us today, almost 100 years after his death. How could you even imagine—or how could you even manufacture—a greater legend?

**BERNDT:** At [Ramanujan's] death, [his wife] Janaki donated all of her husband's papers to the University of Madras in exchange for a small monthly pension for the remainder of her life. I should remark that she told me that during her husband's funeral, some of his papers were stolen. Then in 1923, the library at the University of Madras—actually, it was the registrar to be perfectly correct—sent to Hardy a package of material to be included with the publication of Ramanujan's collected papers. This package included hand-written copies of the three earlier notebooks, and many other papers which were not described in the letter itself, and which were not copied. The registrar, whose name is Francis Dewsbury, asked that these be returned after Hardy made use of them.

As it turned out, when Ramanujan's collected papers were [compiled], neither the notebooks nor any other papers were [included] along with the published papers. So we're quite confident that this shipment of papers to Hardy contained the lost notebook, and probably other papers as well. Hardy evidently had this in his possession for many years, and gave it to G. N. Watson probably sometime between the late 1930s, after Watson's interest in Ramanujan had declined, and Hardy's death in 1947.

**ONO:** G. N. Watson was a distinguished British mathematician who is perhaps most well-known for the books he wrote with E. T. Whittaker. There's a famous book on analysis written by Whittaker and Watson [16]. It's a classic. We still use it today, and Watson is well-known for being one of the co-authors of that book. Watson is also well-known for having been the president of the London Math Society. Years of his life were devoted to studying the mathematics inspired by Ramanujan. Watson was designated to be the first editor of Ramanujan's notebooks after his death [along with B. M. Wilson] [4].

**BERNDT:** So there it was with Watson then for many years—it would have been between 20 and 30 years I guess—before it was found in Watson's attic floor office after he died. And it was found by Robert Rankin of the University of Glasgow, and J. M. Whittaker, the son of E. T. Whittaker who was Watson's co-author.

**ANDREWS:** J. M. Whittaker is a British analyst who was asked by the Royal Society to write the obituary when Watson died. He went to visit Watson's widow and he was shown Watson's study where there was a huge number of papers. He says it was just good luck that he found what I've called the lost notebook, because he only could save a few things and most of it was going to be burned in a few days—he said these papers were about to be burned. Watson kept everything, and so the only way you could deal with stuff was just to trash it.

**ONO:** My guess is that Watson studied these pages, and when he died, he presumably just forgot to point out the importance of these notes. As I understand it, the famous mathematician Robert Rankin saved Watson's papers from incineration. However, we don't know if Rankin had taken the time to study carefully the contents of the estate, the papers.

**BERNDT:** I knew Rankin very well. In fact, after I got my Ph.D. I went to the University of Glasgow for a year, because Rankin was there.

**ONO:** Robert Rankin would be my image of an English gentleman from the 1940s—although he was Scottish. Dressed properly in a coat and tie every time I met him, his posture was perfectly upright, he had the thick accent, he was a proper gentleman. Almost, you know, pipe-smoking, a top hat—that's what we're talking about. Rankin was famous for developing much of the modern theory of modular forms [9].

**BERNDT:** At that time I was strongly interested in modular forms. I might remark that one day I was in Rankin's office—I got my Ph.D. in 1966, so I don't remember if it was 1966 or '67—but anyway, somehow we got talking about Ramanujan and he said, "I have a copy of his notebooks here in my office. I'll be glad to loan them to you if you are interested."

And I said no, I wasn't interested in them. I don't know if it was a mistake or not. Then seven years later, approximately, I did become very strongly interested in Ramanujan. Looking back, I could have learned much more about the lost notebook, because that was the period in which he was actually sorting through the papers left by Watson after he died.

**ANDREWS:** There are many interesting things in [the lost notebook]. Probably the most interesting aspect is that this contains all of his ideas about the mock theta functions. Ramanujan thought these up in the last year of his life and he wrote down many things about them. We've known about them throughout the twentieth century, because three months before [Ramanujan] died he wrote a letter to G. H. Hardy describing them, but not giving nearly the amount of information that was in the lost notebook.

**BERNDT:** Robert Rankin had sorted through all of Watson's material, and sent what we now call the lost notebook to the library at Trinity College on December 26, 1968. And there it sat for I guess about seven-and-a-quarter years before George Andrews re-discovered it.

**ONO:** It's kind of like the trip back in time that none of us could ever take. But if we could have the opportunity to talk to Ramanujan, and talk about his work, the discovery of the lost notebook was something like that—it was like a time machine. Obviously, we couldn't ask Ramanujan any explicit questions, but the uncovering of this notebook was something like that.

**BERNDT:** [Ramanujan's original papers] are still there [at Trinity College], right where they were found by Andrews. He had gone there at the invitation of Lucy Slater, and I think he spent some extra time, so to speak, on this trip to Europe, that he had to use up.

**ANDREWS:** My wife and at the time two small daughters were accompanying me on this trip to Europe. And the reason that I went to Cambridge was because I needed to spend three weeks in Europe, because airline fares were so cheap if you spent three weeks or more, and they were very expensive if you spent less than three weeks. That was why I needed to have other things to do in Europe besides a conference that had only lasted one week. So I proposed to go to look at certain things that I knew were in the Trinity College library.

We had [previously] been [to Trinity College] in 1960 or '61, basically 15 years after World War II, so things were still much less modern than one thinks of now. In particular, I can remember riding on steam trains from Cambridge to London [on the previous trip]. They were actually pulled by steam engines. Much of this was the aspect of recovering from the war. The coinage was still the original British coinage instead of the decimalized version that they have today—shillings and half crowns and whatever—so you would often get pennies with Queen Victoria on them, both Queen Victoria as an old woman and Queen Victoria as a young woman. There would be coins that were more than a 100 years old that you would often just have given to you in change in the stores.

**BERNDT:** There are beautiful campuses [at Trinity College], expansive green lawns, beautiful old buildings—lots of buildings that look like cathedrals, not only the churches but other buildings as well—cold, damp maybe.

**ANDREWS:** There was a conference in Strasbourg, France—that was the first week. Then we went to Paris, and I gave a talk at the University of Paris. Then we went on from there, we spent a day or so at Southampton, and I gave a talk at the University of Southampton. And then we went up to London, because I was working on the collected papers of P. A. MacMahon, and I met MacMahon's great-great-nephew there and got some information from him. Then we went on to Cambridge from there.

So we were staying with Lucy Slater, who was a student of W. N. Bailey, a contemporary of Ramanujan at Cambridge. Bailey was an English mathematician. I do not think he had a Ph.D. because most mathematicians at that time considered getting a Ph.D. as sort of irrelevant to being part of the English college and university scene. He eventually became professor at Bedford College in London. He did a lot of work on the sort of mathematics that interested Ramanujan, and indeed, he devoted much of his professional work to studying things that were either directly of interest to Ramanujan, or parallel to his interests. [Bailey] was a college student at the same time that Ramanujan was in Trinity College, so that he actually physically interacted with Ramanujan. He says—in, I think, a letter of his that I have possession of—that he remembers Ramanujan being in class and when Hardy was lecturing and was uncertain of the value of a definite integral he would ask Ramanujan, and Ramanujan would always tell him the answer.

**BERNDT:** Lucy Slater is really one of the names in [the history of] strong people in *q*-series. She was a student of Bailey, and her book on basic hypergeometric series is a very nice one. She really has done a lot of good work in *q*-series, so the invitation was a very natural one between two experts in *q*-series.

**ANDREWS:** She was the nearest person in Cambridge who was somehow connected personally to people going back to the time of Ramanujan. She was the one who had directed me to these boxes at the Trinity College library.

**BERNDT:** As I understand it, she had told him that he might find the manuscripts that were deposited by the late G. N. Watson interesting to examine.

**ANDREWS:** The architect for [the Wren library] was Christopher Wren, so that's where the name comes from [**5**]. It's part of the Trinity College and you don't really see it as separate from the total effect of Trinity College. When you enter it—I suppose if you are sufficiently impressed by the major people who have been associated with Trinity College such as Isaac Newton, of whom there is a bust in the Wren Library, you are impressed with just the sacredness of it. There were tables in the library and the documents were brought out to me in boxes.

**BERNDT:** And so while doing this, Andrews discovered this sheaf of 138 pages, or sides of pages I guess we should say, in Ramanujan's handwriting.

**ANDREWS:** It turns out that these were the things that had been collected from G. N. Watson's study by Robert Rankin and J. M. Whittaker. Rankin and Whittaker were the people who assisted Watson's widow in contributing these papers to the Trinity College library.

**ONO:** Many of the papers in Watson's estate were forwarded to Trinity College library to be archived, and this included, of course, the lost notebook.

**ANDREWS:** So there was Ramanujan's last letter to Hardy, the actual letter was there, and there was a published version of the letter in his collected works. And I had written my Ph.D. thesis on mock theta functions, which is the subject of this letter. So the letter was there. One recognizes the handwriting both from the letter, and from the photostatic copies of Ramanujan's famous notebooks which were put out by the Tata Institute. So it was easy to recognize the material that was Ramanujan's.

**BERNDT:** He knew immediately that this was from the last year of [Ramanujan's] life because Andrews had written his Ph.D. thesis on mock theta functions and several papers thereafter on the subject. And when he saw all these mock theta function formulas in the manuscript, he knew this was from Ramanujan's last year.

**ANDREWS:** But as I thumbed through these pages, there are very few words in the lost notebook. In particular the phrase "mock theta function" is not in the lost notebook. But the mock theta functions are there, and as soon as you see them, if you know what they are—I was probably the only person on the face of the Earth at the time who knew what they were—Rademacher had died by then, so I was it—…

**ONO:** Rademacher was a German number theorist who spent the bulk of his career as a mathematician at the University of Pennsylvania in Philadelphia. He was trained in the traditional German style of analytic number theory, dating to the ideas of Landau and Siegel. And Rademacher was George Andrews's Ph.D. advisor.

**ANDREWS:** …So once you see [the mock theta functions], you know then that these were things written in the last year of Ramanujan's life, because Ramanujan said in his last letter written three months before his death, "I have discovered some very interesting functions recently which I call 'mock theta functions.'" And he gave examples of them in his letter, and there they were in this manuscript, which meant that this was the stuff he'd been doing in the last year of his life. And for the bulk of the twentieth century, all we knew about what he'd been doing was from this letter that he'd sent to Hardy.

**BERNDT:** Ramanujan did not even give an exact definition of a mock theta function, and it is very hard to give an exact definition. Ken Ono and Kathrin Bringmann in their work [**10**] have a precise definition of a mock theta function using harmonic Maass forms.

**ONO:** I'm very proud of the fact that I've been able to contribute to the story. It's one of my greatest fortunes in my professional career.

**BERNDT:** [Bringmann and Ono's] work was motivated by Sander Zwegers [**18**], but actually Zwegers did not have the definition. Zwegers was a student of Don Zagier at Utrecht in the Netherlands. I was actually on [Zwegers'] thesis committee, and at that time [2002] I myself did not appreciate fully Zwegers' thesis [connecting mock theta functions to harmonic Maass forms]—I mean, I thought it was a great thesis, but I didn't know that it actually would lead to what it has led to.

However, [using harmonic Maass forms] is not how Ramanujan tried to define a mock theta function. Let me just say that they are $q$-series, and they satisfy transformation formulas which are of a similar nature to theta functions. They have some of the behavior of theta functions as $q$ tends to a root of unity.

Generally, $q$-series are infinite series with the following products in the summands: $(1-a)(1-aq)(1-aq^2)(1-aq^3)\ldots(1-aq^n)$, where $n$ is the summation index. Euler was the first [contributor to the theory of $q$-series historically], then L. J. Rogers in the late 19th and beginning of the 20th century, and F. H. Jackson, an English mathematician, and then George Andrews. George Andrews is who I would say is the leading authority on $q$-series, not only now, but in all the history of mathematics—except for Ramanujan.

**ONO:** [Defining $q$-series precisely] is a very hard question, because on the one hand, most mathematicians would agree that any formal power series which can be written in the variable $q$— which any power series could be—should be called a $q$-series. Now, when [Berndt and Andrews] use the terminology "$q$-series," they have a very clear picture of what they have in mind. A Fourier expansion of a modular form is a $q$-series. A hypergeometric series is a $q$-series. A basic hypergeometric series is a $q$-series. The infinite products of the type that Dedekind and Euler would have studied, they are examples of $q$-series. So when people speak of $q$-series, they're really talking about the body of work that's covered by these very specific types of power series. And they are speaking of the combinatorial properties that make up the theory of these objects.

[What $q$ is] is very vague—$q$ is like $x$, it's a variable. Or sometimes $q$ is just a place-holder, that allows us to say the coefficient of $q^n$ just happens to be the $n$th term of this power series. So, often $q$ is just a formal parameter. Sometimes $q$ has to be a complex variable with norm less than or equal to one. And sometimes $q$ has to be $e^{2\pi iz}$, where $z$ is a point in the upper half of the complex plane. There are circumstances where each of those three meanings is the meaning for $q$. You have to be careful when you are reading papers about formal power series or $q$-series, to have a very clear picture of what is meant by $q$ by context.

**BERNDT:** However, [$q$-series] may not have such products [as are mentioned above] in their summands. Often in the theory of $q$-series, one lets the parameter $a$ tend to zero or infinity, and then the products, you might say, disappear—you have to take limits. So you might end up with the series without these products in the summands, and you could get a theta function. Actually what happens is that often, if you can write $q$-series in closed form, the closed form expressions often involve theta functions. This is primarily why theta functions are part of the theory of $q$-series, because of their ubiquitous nature in evaluations and transformation formulas.

**ONO:** A theta function is the prototype of what we now call modular forms, which can be thought of as complex analysis versions of the trigonometric functions like sine and cosine. The theta functions transform back to themselves under what are called Möbius transformations, and these Möbius transformations can be thought of as the complex analytic versions of translations which are central to building trigonometry.

**BERNDT:** The first people who studied [theta functions] were probably Euler and Gauss. Gauss was the first person to actually prove some substantial theorems about them. For example, what's called the Jacobi triple product identity is actually originally due to Gauss, not Jacobi, and that's perhaps the most fundamental property of theta functions.

**ANDREWS:** Jacobi was a major mathematician of the early 19th century, probably the person who really founded, or developed, [the theory of] elliptic functions.

**ONO:** I think of Jacobi as the father of the theory of theta functions.

**ANDREWS:** His *Fundamenta Nova* is the original treatise on theta functions and elliptic functions, and forms the basis of tremendous work on those throughout the 19th century, which expanded into the study of modular forms in the late 19th and then into the 20th century.

**ONO:** I like to view the Jacobi theta function as the first non-trivial example of a modular form. If you go back in time—and you don't have to go very far back in time—when mathematicians discovered the trig functions—sine, cosine—and later discovered their hyperbolic versions, there were immediately great advances in physics. So much of physics is just not at your disposal unless you know the language of trigonometry. Trig functions allow you to describe various types of motion.

Now, in the 19th century mathematicians discovered, for example through solving the heat equation, that there was a need for studying modular forms, these theta functions. And modular forms, just like the trig functions, turn out to be highly symmetric functions which build fields of functions on curves, which we've come to understand are very important in math. These modular forms are functions that satisfy transformation laws—a lot more complicated than what you find for trig functions, but much in the same way.

**BERNDT:** Theta functions converge for absolute value, or modulus, of $q$ less than 1, and one can determine their behavior as $q$ tends to a point on the circle of convergence, say a root of unity. And a similar thing happens for mock theta functions—as you approach a point on the unit circle, the behavior of a mock theta function is similar to a theta function. But as Ramanujan emphasized, you cannot describe the behavior of a mock theta function by just one theta function for all roots of unity. You might mimic a theta function, say for one root of unity or something. So they're similar to theta functions in many respects, but they're much more complex objects, difficult to get your hands on in terms of the properties.

**ONO:** So the question is, what are the mock theta functions? The mock theta functions are pieces, they are pieces of modular forms which do not satisfy the usual laws of calculus coming from complex variables. What we have learned [from the lost notebook], is that Ramanujan defined the mock theta functions—in his last year of life—to be these crazy-looking formal power series, which seem to make no sense in the context of modern complex variables. You couldn't differentiate them and get anything related to a modular form. The reason for that is—we now know—that the mock theta functions are only half of the non-holomorphic modular form. It turns out that classical modular forms are just shadows of harmonic Maass forms.

We use the language like "shadow" all the time in the subject. If you are trying to deduce some information—well, what do they do in astronomy? They're trying to study an object that's many light years away from Earth, a thousand light years from Earth, you want to study some system way out there. How do you study that system? You can't go there. You can only infer properties about that system by measuring the light that passes from that system that eventually reaches us, and we can infer properties about that system by the diffraction pattern that we see when we collect light from these sources far away. That's a shadow.

Hans Maass was a German mathematician, who in the 1930s was the first person to recognize the need for defining these non-holomorphic modular forms, modular forms that are real-analytic but not complex-analytic. It turns out that these harmonic Maass forms are infinite-dimensional spaces of functions that transform like modular forms, because after all, they include all modular forms. In many cases, there are pairs of problems that line up in number theory, where the modular form only has the right to tell you one half of the story, because you're studying the shadow of a harmonic Maass form.

Ramanujan somehow, through some inspiration, was able to discover pieces of functions that people would define in the 1930s. Long before this, in 1920, Ramanujan had discovered about 20 examples of these functions. What is bizarre, and what is befuddling to us, is that he was able to find what most of us would think of as the most complicated part of a harmonic Maass form, without knowing what in the world a shadow was. That's amazing.

**ANDREWS:** The electric moment of realizing that what I had in my hands was the output of Ramanujan in the last year of his life—and nobody had really understood anything about this until the moment I grabbed it—sort of blotted everything else out from my mind.

**ONO:** I can definitely say that had George not found that lost notebook, many of us would be doing different things. I'm pretty sure that maybe a fourth or a fifth of the papers that I've written would not exist, either because of inspiration from the lost notebooks, or following up on the inspiration of others, other people who were inspired by the lost notebooks. Mathematics as we know it—as I know it—would be very different.

**BERNDT:** In the earlier notebooks, out of the more than 3200 results that I counted, there are no more than a half dozen that are incorrect. In the lost notebook, there are more things that are incorrect. There are a number of things which are speculative, so I'm sure that Ramanujan realized that he didn't really have theorems. He was sort of writing down some things—but he really didn't necessarily believe 100 percent that these things were true. He was also dying, and so there were certain formulas that were sort of incomplete. Maybe he scribbled things down and he didn't scribble everything down in the formulas. So it's understandable that there are, you know, mistakes—but I hate to call them "mistakes," because they really just arise from the fact that he was running out of time.

**ONO:** Paper was very expensive in India, and so he did most of his work on a heavy black stone slate. As we understand it, when his calculations resulted in something he wanted to remember—a formula or an expression he wanted to remember—he would then record those findings on the pages which we now refer to as Ramanujan's notebooks. He didn't do his work in a library, he didn't do work at a customary desk. The image should be of Ramanujan sitting on a bench or the stone floor of his temple, his legs [crossed], with slate on his lap madly scribbling with two sheets of paper next to him—not in certain terms of what our image of a notebook is. Our image is as a marbled composition notebook or a spiral notebook. No, no, no. You'd have to think of these as loose sheets of paper which would later be bound.

**ANDREWS:** I was quite excited. I was quite excited. I'd say the main thing it caused me to do was to think about how I could get a copy of this.

**BERNDT:** (Laughs) I remember he told me that he wanted a photocopy made, and the librarian said, "Oh, we can make you a photocopy—it might take about one or two weeks to do that." However, [Andrews] came back the next day and just happened to pass the librarian's desk, and saw that there was a package addressed to him. He said that it took some persuasion on his part to convince the librarian or the clerk at the desk that he could just take this himself and didn't, you know, have to wait for the post man to cart it away and mail it to him. So he could then take a copy of the manuscript back to the U.S. with him.

**ANDREWS:** Probably my wife was the first one I told, and probably next, Lucy Slater, because these were the two people who would be interested in it. After my discovery, I knew that I had really discovered something major.

**ONO:** I think we should make it clear that the lost notebook was never really lost, it was forgotten. It had to be *re*-discovered.

**BERNDT:** The two primary writers of obituaries of Watson, namely Rankin and Whittaker, the two people that visited Mrs. Watson, never did mention the manuscript in their obituaries—so in that sense that knowledge *was* lost, in the fact that it wasn't mentioned in the primary accounts of Watson's papers.

**ONO:** If the contents had been reviewed, the mathematical meaning was never shared with the rest of the mathematical community.

**ANDREWS:** I was so pleased I took my wife and daughters boating, "punting on the Cam." A punt is a long, flat boat and the person directing the punt stands at the rear of the boat. So we went punting on the Cam, and something I forgot about [until recently] was that when you are using this long pole to push the boat along, sometimes it gets stuck in the mud. So it got stuck in the mud and pulled me into the water. [This was] probably the same day, or if not, the very next day—it was a punting trip to celebrate this [discovery].

I was soaking wet and I climbed out of the water. But I dried off enough that we could go. And so, once I was not visibly dripping, we took a cab back to [Slater's] house, and the cab driver was extremely unhappy to see that the seat of his cab was dripping wet, because while I looked that I had dried off, it was only externally. The water from my clothes had covered the seat of his cab.

That evening at Lucy Slater's house, [Slater] took my wife to dinner at her college. The colleges at that time were unisex so that you had to take women to women's colleges and men to men's colleges. She could not take me, so her words to me were, "When a mathematician comes to Cambridge, his host usually takes him off to dinner at the college, and his wife stays home to babysit. Tonight I'm taking your wife to dinner at the college, and you're staying home to babysit. " So she left me there with a magnificent collection of letters from many of the contemporaries of Ramanujan and people of the first part of the 20th century. I definitely had a much better time than I would have going to dinner at any college in the University at that time.

So these are the stories that surround [the discovery of the lost notebook]. That was very close to the end of the trip. We went back to London because we were flying out of Heathrow [Airport]. There was a man that I had corresponded with, Barry Hughes, and I think we stayed at his house, and then went to Heathrow and flew from there.

**BERNDT:** I might have learned of [the discovery] from Richard Askey first, before Andrews mentioned it.

**ANDREWS:** When I got back to Madison, Wisconsin, my host was Richard Askey. Askey is a close friend and collaborator of mine. He spent his [career] at the University of Wisconsin. He is an analyst whose main interests are what are called special functions, which is an aspect of complex analysis. He and I have collaborated over time basically relating the sort of work that he does to the types of things that Ramanujan does.

I had by then a Xerox copy of this. I said [to Askey], "I have a 100-page, unknown manuscript of Ramanujan's in my briefcase, and you can have a look at it for a nickel." (Laughs) I actually don't think I adequately exploited it, though. It was quite a hit.

The first talk I gave on this I think was in Madison, about a year after I had discovered these [notes]. I had worked on them long enough that I had proved enough to really make an interesting talk. I think that was the first [public announcement of the discovery]. And then I gave an invited MAA lecture in Seattle in the summer of 1977, shortly after the Madison talk, so that was really the most widespread talk. D. H. Lehmer's spouse Emma Lehmer, a significant mathematician in her own right, introduced me, saying that the only comparable thing she could think of would be opening an old box and finding Beethoven's Tenth Symphony.

**BERNDT:** Then, in the fall of 1978, George Andrews visited Illinois; in fact, he and his wife stayed at our house, I remember.

**ANDREWS:** [Berndt] and I are collaborating on the published versions of the lost notebook. He is actually my academic great-nephew. His thesis advisor's thesis advisor was also a student of Rademacher.

**ONO:** Rademacher is well-known for many things—he was most well-known for refining the work of Hardy and Ramanujan concerning asymptotic formulas for the partition function. What he did was rather incredible. Instead of just refining the beautiful ideas of Hardy and Ramanujan, he perfected them, and showed how to bootstrap their ideas over and over again to obtain an exact formula for the partition numbers. He described the partition numbers as a convergent infinite sum of crazy transcendental numbers.

**ANDREWS:** So I'm [Berndt's] great-uncle academically.

**BERNDT:** [Andrews] told me that Watson and Wilson's efforts to edit the notebooks were actually preserved at the library at Trinity College at Cambridge. He had discovered this on his visit there. So I wrote for a copy of all these notes that Watson and Wilson had accumulated, and I thought, well maybe I could just go back and edit further chapters.

As I was sort of winding up my work [editing] the [earlier] notebooks, I also started to work on the lost notebook, I would say mid-1990s, although it might have actually been even earlier—it might have been about '93, I'm not sure. If I was stuck on something I was doing in the earlier notebooks, then I would start working on the lost notebook.

**ANDREWS:** My own work, really since I was a graduate student, has been more or less in the area of things that Ramanujan would have taken an interest in—the theory of partitions, these $q$-series, and so on—all things that were grist for Ramanujan's mill.

**BERNDT:** There are a lot of things in the lost notebook that weren't at the forefront of George Andrews's interest, that were at the forefront of mine, because they had connected with things in the earlier notebooks [**7**]. In particular, things involving modular equations and theta functions and stuff involving analytic number theory and analysis—that is, classical analytic number theory and classical analysis.

**ANDREWS:** Analytic number theory is, succinctly, the utilization of complex analysis to find out facts about numbers.

**BERNDT:** I had seen these [papers] some years ago, and I was really anxious to get to work on them. But I didn't really want to start full-time because I really wanted to finish the earlier notebooks. But I did sort of start a little bit, and once the earlier notebooks were done, then I went to work full-time. In other words, there was really sort of never any argument between us about what we were going to work on. And it was clear that for a lot of these things George is the leading expert in the world—but on the other hand, there were certain things that he was not so interested in, and these happened to be things that I was interested in.

**ANDREWS:** It is amazing, a great opportunity for those of us who are less brilliant than Ramanujan. I'm grateful for all of it.

**BERNDT:** It's just been a privilege and honor for me to work on these formulas of Ramanujan. And I'm very fortunate to collaborate with George Andrews on this. It's been a privilege and very enjoyable for me. It was fortunate for me there was a lot left to do after Andrews worked on the lost notebook for 20 years—there was still a lot left that no one had examined.

**ONO:** We're still busy. We've got a lot to do.

**BERNDT:** It's been said that the discovery of the lost notebook is comparable to discovering Beethoven's Tenth Symphony, and I think that certainly is correct. It was a startling event, and really a lot of interest certainly arose from this, but I think the interest has really increased over the years.

Because, you know, first when you hear about such a discovery, you say, "Oh, yeah, that's wonderful, that's great!" But then, as you learn what's in it, and then as people learn, "Oh, this is connected with my work," then the interest actually increases and balloons. I think that's actually what's happened. You know, especially in the last half dozen years, a lot of people have become more interested in the lost notebook when they didn't have any particular interest, when they see, "Oh, this is connected with things that I'm interested in!"

**ANDREWS:** Mathematically, these things have led to really exciting and interesting discoveries, so it was something very, very fortunate from every point of view.

**ONO:** There are a lot of things that would have never happened. People who are studying the calculation of black hole entropies, the people at Harvard who are studying moonshine for the M24 group—none of that could have happened. People who are studying differential topology of the complex projective plane wouldn't be able to even begin to state some of the formulas that are now central to their work. Now that we know the mathematical meaning of the term "mock theta function," we've been able to develop general theory, and over the last five or six years we've begun to understand that this theory has tentacles that reach out to many different areas of math and science.

**BERNDT:** This excitement is still at a high pitch, a high fever. The interest has been growing and growing as we learned more about what's in [the lost notebook].

**ONO:** People quickly began to prove theorems in number theory, like the theory of partitions, by using what we've now learned about the mock theta functions. But these other areas of science in which the mock theta functions play an important role, they make up quite a long list. They include things like probability theory, people that study the thresholds of certain cellular automata by using asymptotic properties of mock theta functions—a cellular automaton can be viewed as being something like a toy model for a cancer tumor, a probabilistic model that varies over time, and one wants to decide what the initial parameters of such a model should be before we have confidence that a space concerned can be dominated by a tumor or not.

[Other] applications include geometric topology, differential topology, where we're trying to understand whether certain objects like "two-form" four-manifolds can be continuously deformed into each other—this is a very difficult problem. Ten years ago nobody would have thought the mock theta functions would have anything to do with understanding objects in space and their interrelationships.

Moonshine is a statement about understanding the representation theory of simple sporadic groups. They are very bizarre [groups], and they are realized in nature in ways which just don't seem to fit in with the other groups. They've stood out as being objects that are clearly somehow central to mathematics. We need to understand every facet of these groups if we really want to be in a position where we can say that group theory is, in some sense, solved—and if we really want to be in a position where we can apply group theory in its full glory.

And this is where moonshine comes in. Moonshine is a procedure by which we can relate all of the irreducible representations of a given finite simple group to a power series. People like John Conway and others have done this. These moonshine functions also appear in physics, and because of that people are quite interested in them. It turns out that some of the mock theta functions also appear to be moonshine functions. So people like [Amanda] Folsom and Bringmann are studying this kind of moonshine, and it's quite striking how the coefficients of the mock theta functions appear to be related to representation theory. That certainly comes as a surprise. Nobody would have guessed that 20, 30 years ago, or even ten years ago.

It's also in this connection that the entropy of black holes comes into play. I'm not an expert on this, but the entropy of black holes is directly related to the crazy infinite sums that Rademacher invented to write down his exact formula for [the partition function]. Formulas like this now appear in the calculations of black hole entropies. These Rademacher sums, when they are related to the mock theta functions, seem to correspond to some of the sums that people have been producing from the physics.

**BERNDT:** Just in the last couple of years I've made some discoveries that are really remarkable— things that I didn't know were in the lost notebook. The most recent [discovery] was just this past fall. I'd seen this two-page manuscript [of Ramanujan's] and I'd passed over it many times. But now I'm trying to get everything done in the lost notebook, so I said, I have to look this over carefully and see what this is. And I discovered that Ramanujan had found the best Diophantine approximation to the exponential function. He wrote it as $e^{2/a}$, where $a$ is any nonzero integer. In particular, this approximation for the exponential function, this Diophantine approximation, beats the classical Diophantine approximation.

**ANDREWS:** Diophantus was a Greek who was probably viewed as the father of number theory, in that he looked at solving various types of equations in integers. And Diophantine approximation is approximating numbers—algebraic numbers, transcendental numbers—by rational expressions.

**BERNDT:** Then I found that the best Diophantine approximation had not been published until 1978. This was done by an Australian number theorist by the name of C. S. Davis. So here Ramanujan had, somewhere between 1915 and 1920, the best Diophantine approximation of the exponential function, and it was not discovered by others until 1978. Ramanujan not only gets the Diophantine approximation, but proves that it is the best that you can get. Now, to tell you the truth, the proof was actually wrong, but I was able to fix it without much difficulty. There are three of these partial manuscripts in Diophantine approximations.

What happened [next] is that during the fall, I went to my colleague Alexandru Zaharescu and I said, "Look what I discovered—this wonderful Diophantine approximation of the exponential function of Ramanujan!" At that moment I didn't actually know it was the best—well, I knew it was the best because Ramanujan had proved it, but I didn't know about Davis' paper at that time. Zaharescu said, "Oh, that's exactly what I need! I'm working on this problem, this open problem of Jonathan Sondow." So this was exactly what Zaharescu needed in his work on this problem.

I discovered something else [in the lost notebook] recently that I had previously discussed with Andrews. He had said, "Oh, that's an integral involving mock theta functions." So I've been looking at it lately, and actually it's not involving mock theta functions. It involves maybe something else which is new, that I have no idea what this is.

**ANDREWS:** It's very exciting. This is very exciting, and it's a tribute to the genius of Ramanujan that the things that he thought of in this last year of his life when he lay dying, are now becoming really major hot topics, 90 years later.

**BERNDT:** A lot of times, when people ask mathematicians what they're doing, and how they discovered it—it's very hard to actually convey these ideas to lay-people, or even to other mathematicians.

**ANDREWS:** What stick out in my mind, that have come up time after time, are discussions about [Ramanujan's] conversations with his friends, his telling them that the goddess Namagiri whispered formulas to him in his sleep, and when he woke up he would write them down.

**BERNDT:** Ramanujan could say, "I got these ideas from the goddess Namagiri in dreams," and we might say that, you know, I got these from jogging, or hiking up a mountain trail, or taking a shower, or something like this. In other words, "This is my inspiration." You know, sometimes it just hits you after you've been thinking about things for a long time.

**ONO:** Ramanujan as an untrained mathematician, as an amateur, didn't recognize the need to provide rigorous proofs of his findings, because his findings came to him as visions. For us as Western mathematicians, we have the "aha," light bulb moment when you finally recognize or understand or know how to solve a problem you've been struggling with. For us, we come to view that as finally realizing the fruits of hard effort, when you finally understand a problem. I'm confident that when Ramanujan had his "aha" moment, that he viewed this as a gift, a gift from the family goddess, the goddess Namagiri of Namakkal. I'm sure that description is accurate.

**ANDREWS:** Hardy did not believe that Ramanujan had any serious religion, but was only an observant Hindu. Hardy bases this on the fact that once at tea, Ramanujan said to Hardy that he believed all religions were more or less equally true. Of course Hardy, with his solid Western logic, immediately recognized that, since if two religions contradict each other and they are both equally true, that means they're both false. And therefore what he viewed Ramanujan as saying was that all religions are false, or not to be taken seriously.

On the other hand, if you think of Ramanujan as a young man at tea with his mentor and benefactor G. H. Hardy, who was well known as a vigorous atheist, what would you say when the question of religion came up that would be truthful, and that would put Hardy off without making him think that you were a religious fanatic? There could not be a better answer than, "I regard all religions as more or less equally true." It worked perfectly. So I think he was somebody who was sharp, and a religious Hindu.

**BERNDT:** I strongly believe, as Hardy believed, that mathematicians when you come down to it all think the same. In other words, you do think deductively and you often make these deductions. It just hits you, but often you've thought about these things for a long time before you make the right connections in the neurons up in your brain.

I think Ramanujan was a mathematician like any other mathematician, except that, of course, he had more creativity, and he was smarter than the vast majority of us. But I think people make too much of mystical aspects or, you know, religious aspects of him. I'm not discounting religion at all—I mean, I myself am quite religious—but I think a Hindu doing mathematics is really no different than a Christian or a Muslim doing mathematics.

**ANDREWS:** My feeling is that Ramanujan is genuinely unique in any real sense, so that he is not to be seen as sort of a natural development of anything. Certainly his mathematics develops from a few unusual sources related to Western mathematics, but nothing that would at all suggest that he would have made the discoveries that he made. Many times when I am speaking to an Indian audience, during the question-and-answer period somebody will have an explanation for Ramanujan: it was the vegetarian diet, it was the Hindu religion, it was the Vedas, it was something else to do with the culture of India.

And my response is always the same: the first thing we must say is, we must pay our respects to India. You are way ahead of the rest of the world in producing Ramanujans. You've produced one, the rest of us have produced zero. However, the number you've produced is not enough to draw any scientific conclusions about how Ramanujan was produced, or what explains his magical genius.

**ONO:** Kumbakonam has always been a very significant town to Hindu culture. Every twelve years, millions of Hindus descend on this town of Kumbakonam for what is called the Mahamaham festival. And if you ever visit Kumbakonam, what you'll quickly discover is that this is a town of many, many, many temples—giant, beautiful temples that were built out of stones that were brought to South India by elephants centuries ago. Many of these temples tower hundreds of feet into the sky, they're brilliantly colored—beautiful pastel colors, reds and blues and yellows and orange—where carved on every square inch of the temple you find bits and pieces of ancient Hindu legend, just scrawled and carved pictorially on the sides of these giant buildings. This was the environment in which Ramanujan grew up and lived his life.

This is subject to debate, but for me there is no debate. Ramanujan grew up in South India, in the state of Tamil Nadu, where at the time and even today, well over ninety percent of the citizens are Hindu. I've been there many times. I've been to his temple, I see at the end of the day hundreds of people lining up to participate in their daily ceremonies. Everything about life in South India today, and almost certainly during Ramanujan's day, revolves around their Hindu beliefs. So it's hard to imagine that, as a devout Hindu, Ramanujan would be any other way, that he would be anything other than your image of a devout Hindu from the early twentieth century. What would that mean? He probably participated in these ceremonies at least once, perhaps two or three times a day. We understand that he did much of his work, certainly during the summer months, on the cool stone floors inside the bowels of this giant temple that was one of the cool places that you could find.

I think it would have been inescapable—I think everything about his life would have revolved around his religious beliefs. In fact, every morning [he and his family] would make these beautiful geometric designs out of rice flour as an offering to the family god. Why would they do this? They would make these beautiful geometric, flowery designs on the sidewalk in front of the home out of rice flour every day as a religious offering. You wouldn't do that unless you were really devout. Now, I can't say how he lived while he was in England, because, of course, that would represent a huge culture shock—not to mention that as a vegetarian, you would probably have a very difficult time finding suitable things to eat—but surely Ramanujan was a very devout Hindu.

**BERNDT:** I actually asked Mrs. Ramanujan—Janaki—this question. I said, "We've got different opinions as to whether your husband was religious or not. We've got Hardy and [Indian physicist] Chandrasekhar on one side, and then his Indian biographers on the other side." So I asked Mrs.

Ramanujan this question, and she said, "No he was not religious at all. He never went to the temple, he was too busy doing his sums." That's how Janaki described Ramanujan, he was always "doing his sums." And she said, "Often I had to feed him and put food in his hand so he could continue doing his sums and not have interruptions."

**ONO:** Ramanujan shouldn't be thought of as someone who shed light on deep conjectures, someone who will have theories named after him. That wasn't Ramanujan. Ramanujan was the anticipator. He was the person that could discover beautiful examples, prototypes of general theory that would become the basis of modern number theory. That's what I think of as the lost notebook.

**BERNDT:** These many wonderful formulas of Ramanujan—they're just so beautiful and surprising, and things you would never think of. You constantly feel that if Ramanujan hadn't discovered these things, no one would ever have discovered them.

**ONO:** We don't make mathematicians that way anymore. We're not able to make mathematicians who can anticipate the future. We train mathematicians by teaching them a body of knowledge, exposing arguments that tend to work on certain types of problems. We teach bodies of knowledge, but you can't teach someone to anticipate what the future will be—that's a gift. And Ramanujan was world-class at that.

ACKNOWLEDGMENTS

The author would like to express heartfelt thanks to Professors Andrews, Berndt and Ono for their graciousness in participating in these lengthy interviews, the fascinating contents of which extend many times beyond what is included here. In addition, I would like to thank Benjamin Braun at the University of Kentucky, whose History of Mathematics course provided the initial motivation for researching and writing this article; Neil Calkin at Clemson University for background information about Trinity College; as well as William Dunham at Muhlenberg College, Andrew Granville at the University of Montreal, and my friend, author Benjamin Phelan, for editorial revisions.

DEPT. OF MATHEMATICS AND COMPUTER SCIENCE, EMORY UNIVERSITY, ATLANTA, GEORGIA 30322
*E-mail address:* `robert.schneider@emory.edu`